\documentclass[11pt,a4paper,reqno]{amsart}
\usepackage[utf8]{inputenc}
\usepackage{amsmath, amsthm}
\usepackage{amsfonts}
\usepackage{amssymb}
\usepackage{graphicx}
\usepackage[left=2.5cm,right=2.5cm,top=2cm,bottom=2cm]{geometry}
\usepackage{amscd}
\usepackage[all]{xy}
\usepackage{hyperref}
\usepackage{gensymb}
\usepackage{enumerate}
\usepackage{csquotes}
\usepackage{comment}
\usepackage{hhline}
\usepackage{epstopdf}
\usepackage{float}

\newtheorem{theorem}{Theorem}

\begin{document}

\title{Numerical Methods and Closed Orbits in the Kepler-Heisenberg Problem}

\author{Victor Dods}
\email{victor.dods@gmail.com}

\author{Corey Shanbrom}
\address{California State University, Sacramento, 6000 J St., Sacramento, CA 95819, USA}
\email{corey.shanbrom@csus.edu}
\thanks{This material is based upon work supported by the National Science Foundation under Grant No. DMS-1440140 while both authors were in residence at the Mathematical Sciences Research Institute in Berkeley, California, during Summer 2017.  CS was supported by a Simons Travel Grant from the American Mathematical Society.}

\date{\today}

\subjclass[2010]{65P10, 70H12, 70F05, 53C17}


\begin{abstract}
The Kepler-Heisenberg problem is that of determining the motion of a planet around a sun in the sub-Riemannian Heisenberg group
.  The sub-Riemannian Hamiltonian provides the kinetic energy, and the gravitational potential is given by the fundamental solution to the sub-Laplacian.  This system is known to admit closed orbits, which all lie within a fundamental integrable subsystem.  Here, we develop a computer program which finds these closed orbits using Monte Carlo optimization with a shooting method, and applying a recently developed symplectic integrator for nonseparable Hamiltonians. 
Our main result is the discovery of a family of flower-like periodic orbits
with previously unknown symmetry types.  We encode these symmetry types as rational numbers and provide evidence that these periodic orbits densely populate a one-dimensional set of initial conditions parametrized by the orbit's angular momentum.  
We 
provide links to all code developed.
\end{abstract}

\maketitle


\section{Introduction}

In geometric mechanics one usually constructs a dynamical system on the cotangent bundle of a Riemannian manifold $(M, g)$ by taking a Hamiltonian of the form $H=K+U$, where the kinetic energy $K$ is determined by the metric $g$ and the potential energy $U$ is chosen to represent a particular physical system.  In particular, the classical Kepler problem has been extensively studied in spaces of constant curvature; see \cite{Diacu} for a thorough history.  In \cite{MS}, we first posed the Kepler problem on the Heisenberg group in the following manner.

Let $(\mathcal H, D, \langle \cdot, \cdot\rangle)$ denote the sub-Riemannian geometry of the \emph{Heisenberg group}:
\begin{itemize}
\item $\mathcal H$ is diffeomorphic to $\mathbb R^3$ with usual global coordinates $(x, y, z)$
\item $D$ is the plane field distribution spanned by the vector fields $X := \frac{\partial}{\partial x} -\frac{1}{2}y\frac{\partial}{\partial z}$ and $Y := \frac{\partial}{\partial y} +\frac{1}{2}x\frac{\partial}{\partial z}$ 
\item $\langle \cdot, \cdot\rangle$ is the inner product on $D$ which makes $X$ and $Y$ orthonormal; that is, $ds^2 = (dx^2+dy^2)|_D$.
\end{itemize} 
See \cite{Tour} for a detailed description of this geometry.

We define the Kepler problem on the Heisenberg group to be the dynamical system on $T^*\mathcal H=(x, y, z, p_x, p_y, p_z)$ with Hamiltonian 
$$H=\underbrace{\tfrac{1}{2}((p_x-\tfrac{1}{2}yp_z)^2+(p_y+\tfrac{1}{2}xp_z)^2)}_K\underbrace{-\frac{1}{8\pi\sqrt{(x^2+y^2)^2+{16}z^2}}}_{U}.$$
The kinetic energy is $K=\frac{1}{2}(P_X^2 +P_Y^2)$, where $P_X= p_x-\frac{1}{2}yp_z$ and $P_Y= p_y+\frac{1}{2}xp_z$ are the dual momenta to the vector fields $X$ and $Y$; the flow of $K$ gives the geodesics in $\mathcal H$.  
The potential energy $U$ is the fundamental solution\footnote{Thanks to Michael VanValkenburgh (Sac State) for correcting the value $\alpha=1/8\pi$, which incorrectly appeared as $2/\pi$ in a previous paper; no prior results are affected by this change.} to the Heisenberg sub-Laplacian; see \cite{Folland}.  We use the notation $q=(x, y, z)$ so a solution to this system can be expressed as $q(t)$.

In \cite{MS} we analyzed many properties of the resulting system and proved that the dynamics are integrable on the invariant hypersurface $\{H=0\}$. We also showed that any closed orbits must lie on this hypersurface.
The Kepler-Hesienberg system is at least \emph{partially integrable}, in that both the total energy $H$ and the angular momentum $p_{\theta}=xp_y-yp_x$ are conserved.  Moreover, the quantity $J=xp_x + yp_y +2zp_z$  generates the Carnot group dilations in $T^*\mathcal H$, which are given by $$\delta_{\lambda}(x, y, z, p_x, p_y, p_z)=(\lambda x, \lambda y, \lambda^2 z, \lambda^{-1}p_x, \lambda^{-1}p_y, \lambda^{-2}p_z)$$ for $\lambda>0$.
It always satisfies $\dot J =2H$ and is therefore conserved in closed orbits, for which $H=0$.  

In \cite{CS} we proved that closed orbits exist.
\begin{theorem}[\cite{CS}]\label{periodic}
For any odd integer $k \geq 3$, there exists a periodic orbit with $k$-fold rotational symmetry about the $z$-axis.
\end{theorem}
While crude numerical approximations were given in \cite{MS} and \cite{CS}, these were literally found by guessing, and the images were unsatisfactory as the orbits did not actually close up. 
Moreover, Theorem \ref{periodic} only provides existence -- there was no way to find these orbits in order to further analyze their properties.  Even worse, Theorem \ref{periodic} was only proved for odd $k$.  The original goal of the present work was to address these two inadequacies.  Here, we develop a Python program for finding numerical approximations to periodic orbits enjoying $k$-fold rotational symmetry for any $k> 1$, even or odd.  We also discover more intricate symmetry details; there are in fact $\varphi(k)$ many symmetry classes for any fixed $k$, where $\varphi$ denotes Euler's totient function. See Section \ref{results}.

\textbf{Main Numerical Discovery.} Let $j$ and $k$ be positive integers, and let $\omega=\exp (2\pi i/k)$ be the generator of the cyclic group $\mathbb{Z}_k$, acting on $\mathbb{R}^3$ by rotation about the $z$-axis by $2\pi/k$ radians.  Then we have numerically discovered a large number of solutions to the equation 
\begin{equation} \label{1}
q(t+T/k)=\omega^j q(t).
\end{equation}  

\begin{figure}
\centering
\begin{tabular}{cc}
\includegraphics[width=.26\textwidth]{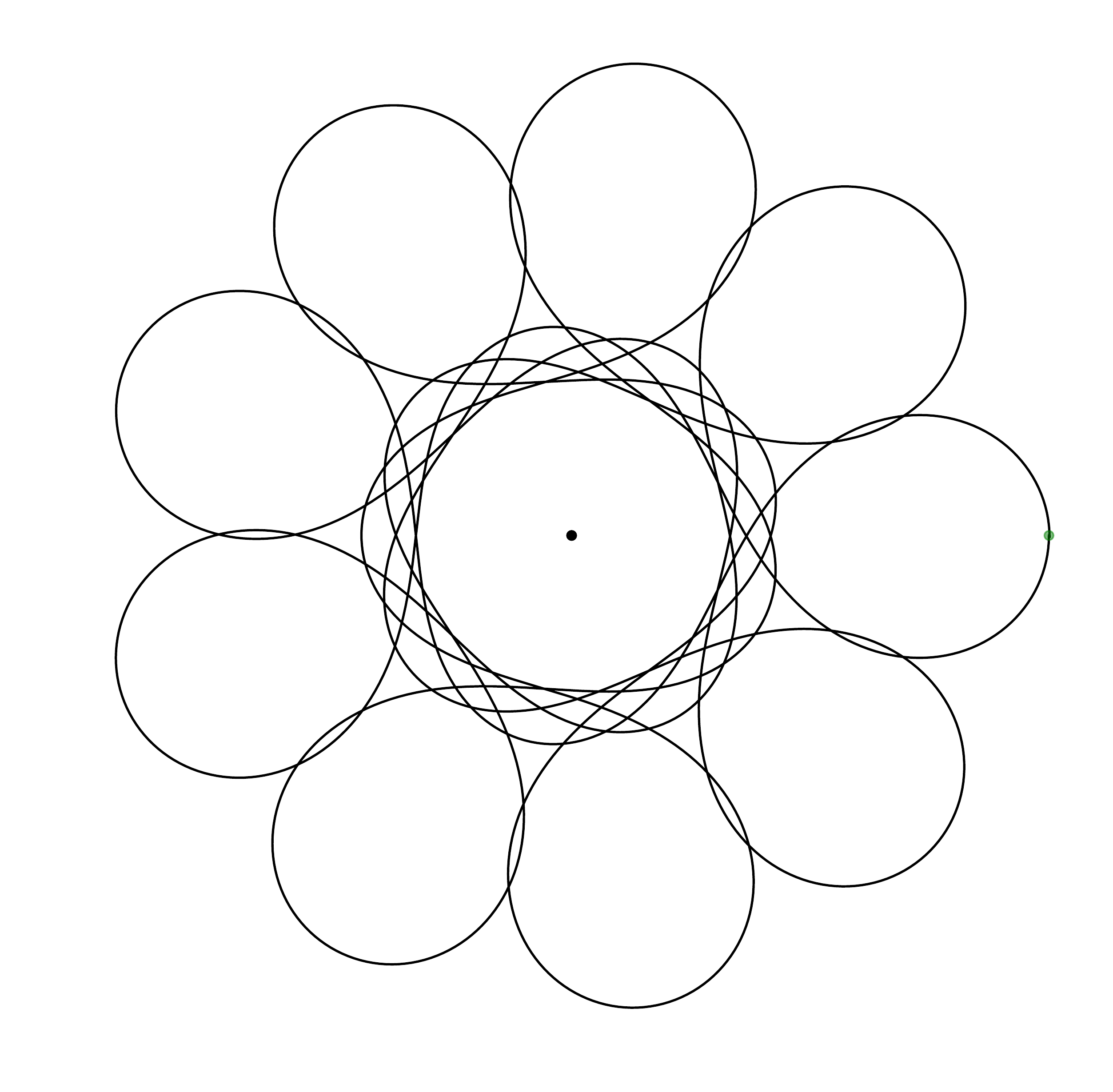}
\includegraphics[width=.26\textwidth]{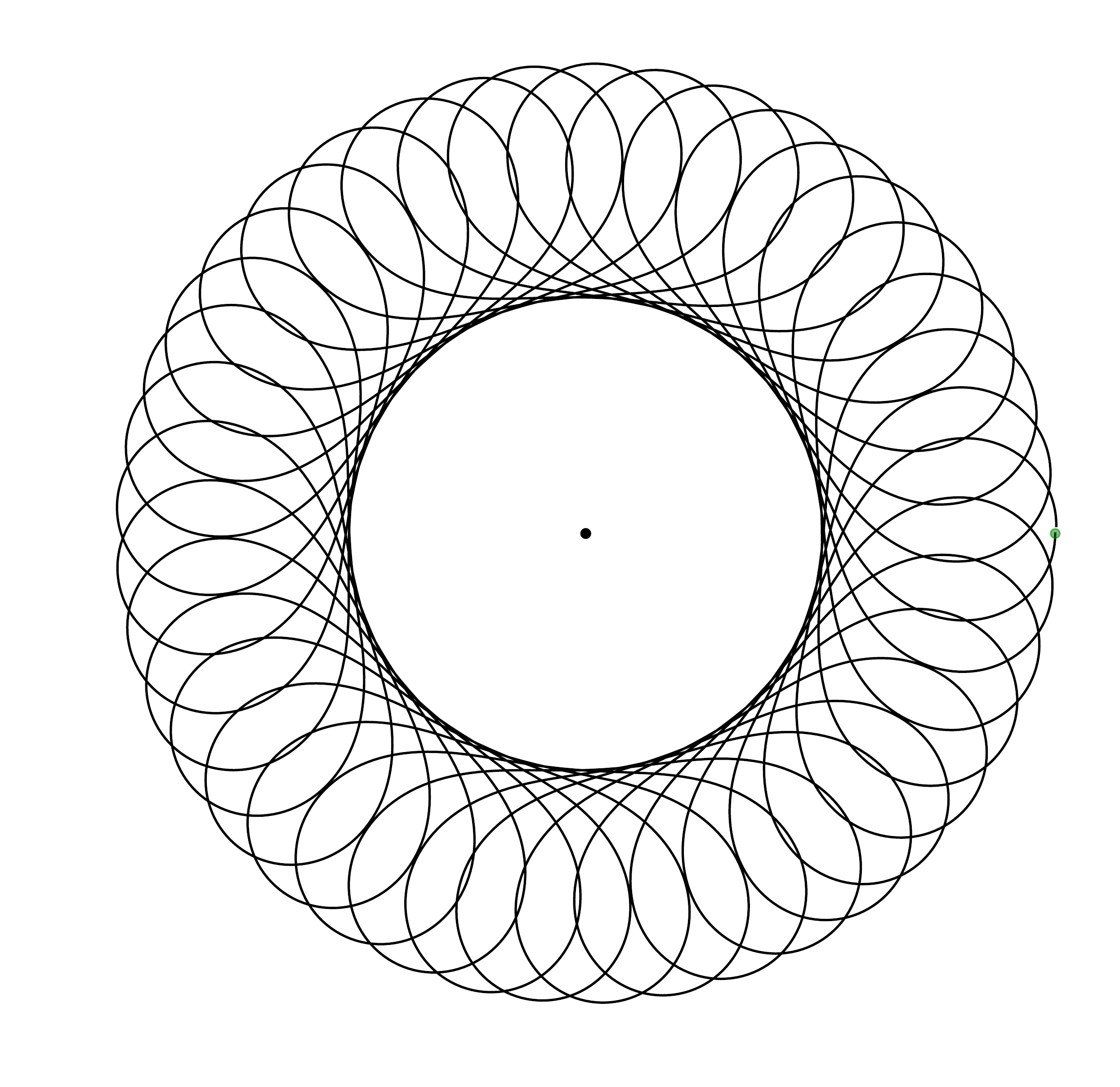}
\\
\includegraphics[width=.26\textwidth]{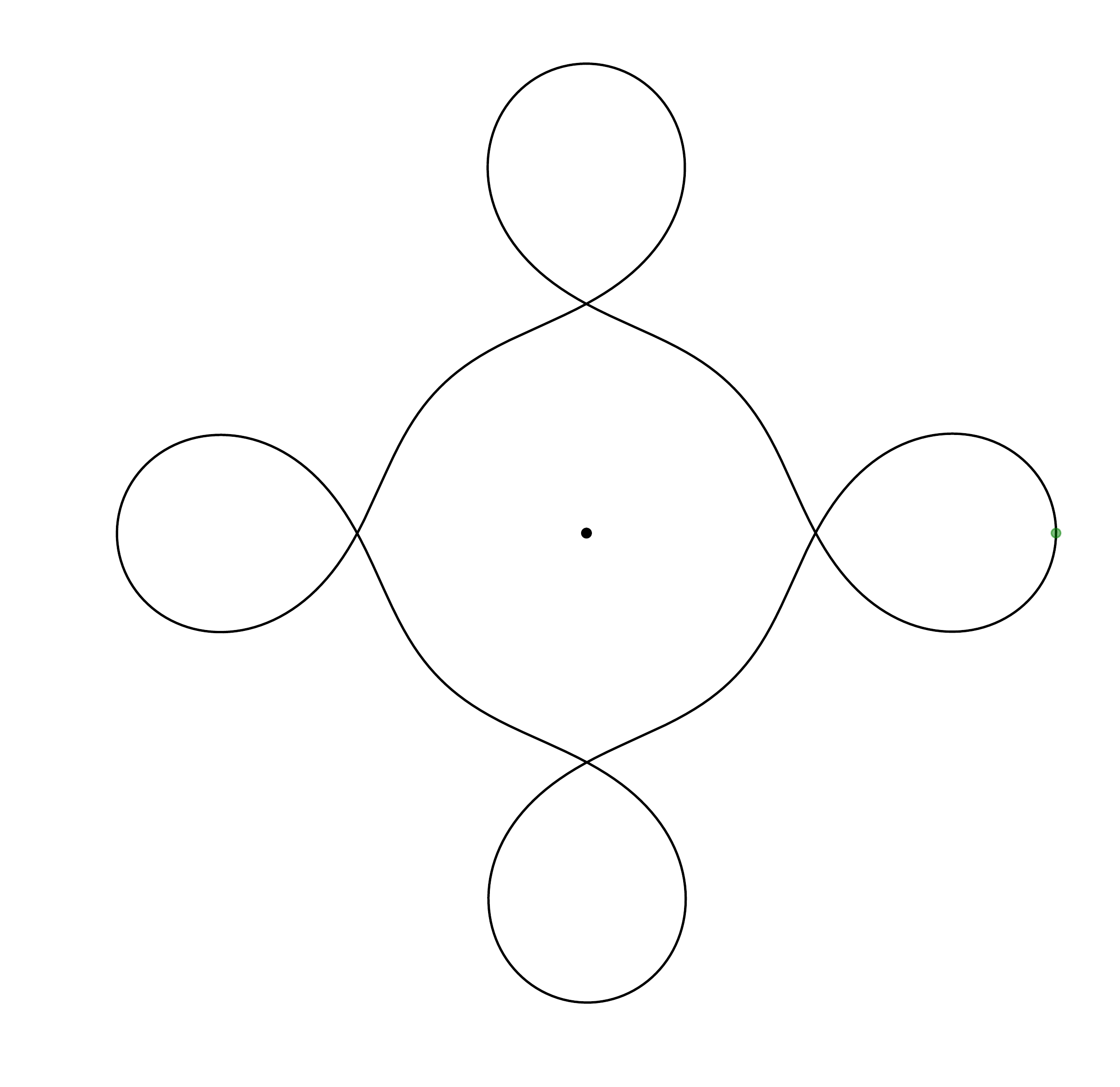}
\includegraphics[width=.26\textwidth]{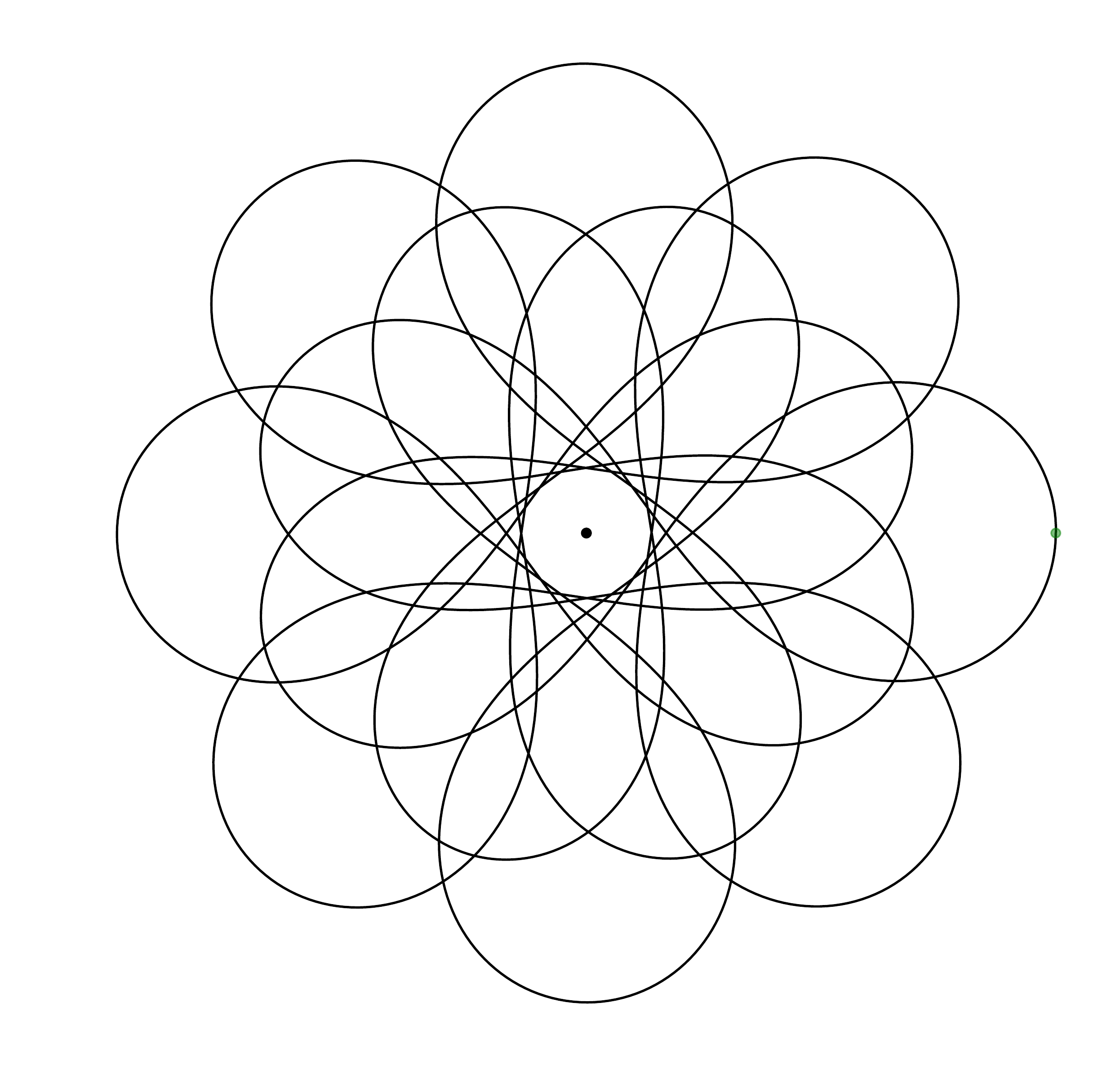}
\end{tabular}
\caption{Closed orbits with $k$-fold rotational symmetry for $k=9, 41, 4, 8$.  In these figures and those following, the black dot represents the origin 
and all orbits are projected to the $xy$-plane.  The $z$-coordinate is equal to the area traced by this projection
.} 
\label{k-fold}
\end{figure}

In \cite{CS}, we proved this result in the case $j=1$ and $k$ odd, with essentially no pictures of the corresponding solutions.
Figure \ref{k-fold} shows orbits with $k=9, 41, 4, 8$ and $j=5, 6, 1, 7$, respectively.  They were all found using random initial conditions, then the optimization procedure described in Section \ref{program}.

\section{Implementation}

\begin{figure} 
\centering
\includegraphics[width=.9\textwidth]{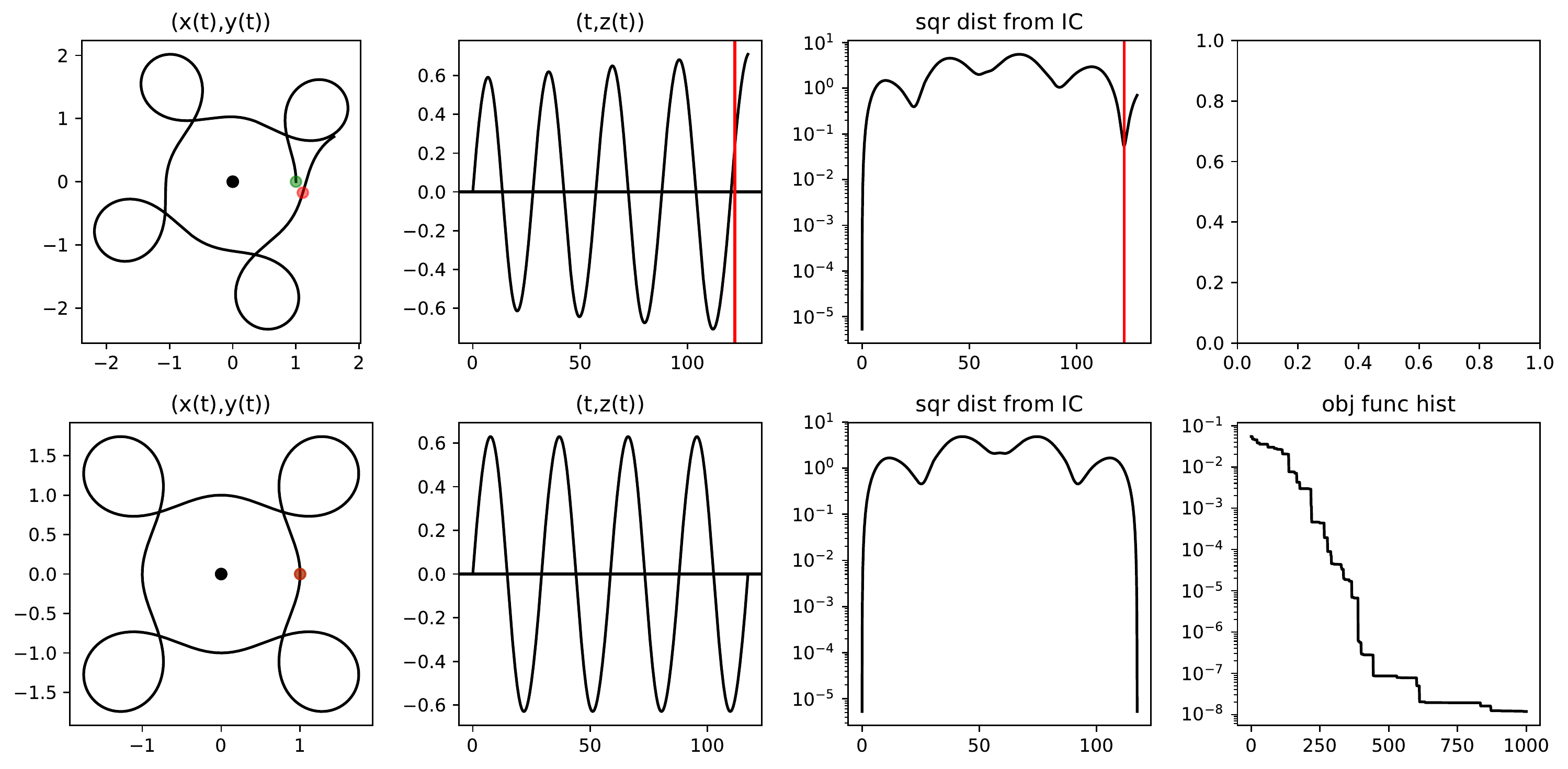}
\caption{This orbit was found using our \texttt{search} procedure (Section \ref{program}).  The top row represents a quasi-periodic orbit with random initial conditions, the bottom row is its periodic neighbor. The data in the columns are described in (\ref{data}).  The bottom right plot shows the objective function decreasing over 1000 iterations of step (\ref{update}).} 
\label{search}
\end{figure}

We use the recently developed symplectic integrator for nonseparable Hamiltonians, as described in \cite{Tao}.  This builds on an innovative algorithm developed in \cite{Pihajoki}, which was first to consider a symplectic leapfrog method for nonseparable Hamiltonians by mixing coordinates in extended phase space.  The method in \cite{Tao} imposes an additional constraint which binds the two copies of phase space, resulting in improved long time performance without loss of symplecticity.  We combine this integrator with our own shooting method for finding closed orbits.  See Figure \ref{search} for an example of a successful search.

\subsection{Program description}\label{program}

Roughly, the purpose of our program is to input a choice of initial condition and  output a \emph{closed} orbit with nearby initial conditions.  Note that:
\begin{itemize}
\item The original choice of initial condition, which we call the \emph{starting initial condition}, must satisfy $H=0$.
\item The symmetries of our system allow the user to simply enter two (of six) coordinates, $p_x$ and $p_y$, if she so chooses.  See Section \ref{reduction}.
\item Alternatively, one can generate random initial conditions (with $H=0$) by simply inputting an integer seed for our (psuedo-)random number generator.  This is our \texttt{search} feature. 
\end{itemize}

Our program follows four main steps.  
\begin{enumerate}

\item Input starting initial condition $X_0=(x(0), y(0), z(0), p_x(0), p_y(0), p_z(0))$.  As noted above, there are three options here: choose all six, choose just $p_x$ and $p_y$, or choose randomly.

\item When using \texttt{search}, the program first determines whether the random starting initial condition is close enough to that of a closed orbit (if the starting initial condition is chosen manually, this step is skipped).  We define the \emph{objective function}, $obj$, for our optimization problem to be the smallest local minimum of the (Euclidean) distance squared from $X_0$ in phase space; this choice makes our approach a shooting method.  Setting a threshold of $obj < 0.1$ effectively defines a small neighborhood around $X_0$.    
We employ the symplectic integrator in \cite{Tao} to integrate the random starting initial condition for some small time. If the resulting curve leaves and returns to this neighborhood, this $X_0$ is satisfactory and we proceed to step (\ref{update}).  Otherwise, we increase the time and repeat.  If the orbit never returns to the neighborhood within some set time threshold then we abort and save the curve in a directory of failures called \texttt{abortive}.

\item \label{update} Once a starting initial condition $X_0$ has been deemed satisfactory, we use a Monte Carlo method to optimize the objective function above.  Use a radially symmetric distribution defined on an annular neighborhood centered at $X_0$ (we used an outer radius of $10^{-1}$ and an inner radius of $10^{-12}$).  Randomly pick a point in this distribution as a new initial condition.  Integrate the dynamics using the symplectic integrator (for the same time as for $X_0$) and evaluate the objective function along this curve.  If this new evaluation is not less than the previous value, repeat with a new random point from the distribution.  If the objective function did improve, repeat step (\ref{update}) with this new initial condition; this is the \emph{update step}.  Iterate this procedure until a specified number of points have been tested (we used 100--1000 iterations).

\item \label{data} The output is an array of figures, where the user can specify which of the following curves are to be plotted for both the starting curve and the optimized curve:
\begin{enumerate}[i.]
\item The curve in configuration space projected to the $xy$-plane.  This is the best representation of the periodicity and symmetry of the orbit.  
\item The $z$-coordinate over time.  This is always qualitatively sinusoidal and can be recovered from the projected curve according to the sub-Riemannian structure.
\item The objective function over time, which has a minimum where the orbit closes.
\item The objective function with respect to the number of optimization method iterations.
\end{enumerate}

\end{enumerate}

Other options are available to the user as well, including plots of $H$ and $J$ over time.  The program also outputs a binary file with all generated data as well as a human-readable summary file. And we developed both random and uniform sampling to search within specified regions of phase space; this was used to generate the images in Figure 
\ref{scatter1}.  Details can be found in \cite{Github}.


\subsection{Reduction of sample space} \label{reduction}

Ostensibly, our initial conditions could lie anywhere in $T^*\mathcal H \cong \mathbb R^6$.  However, we can narrow our search to an embedded 2-dimensional submanifold as follows.  We use basic results about the system, which appear in \cite{MS}.
\begin{enumerate}[1.]
\item By conservation of the angular momentum $p_{\theta}=xp_y-yp_x$, we can assume $y=0$.  

\item All orbits must satisfy the sub-Riemannian (horizontal) constraint $\dot z = \frac{1}{2}(x\dot y - y \dot x)$.  Integrating this condition and applying Green's theorem forces the $z$-coordinate to equal the area traced out by the projection of the orbit to the $xy$-plane. 
Since we are looking for closed orbits, the $z$-coordinate must be zero at some time, which we can take to be the initial time; we thus impose $z=0$.

\item Recall that $J=xp_x + yp_y +2zp_z$ is conserved in closed orbits.
Utilizing this dilational symmetry (which does not exist in any homogeneous Riemannian manifolds besides Euclidean spaces
), we can set $x=1$.  Note that a choice of $x=0$ is impermissible since $(0, 0, 0)$ represents the sun, where our potential energy is singular.

\item Finally, since closed orbits have zero energy, we can solve $H=0$ for $p_z$, obtaining $p_z = -2p_y \pm \sqrt{\frac{1-4\pi p_x^2}{\pi}}$ , which has a real solution if $|p_x| \leq \frac{1}{2\sqrt \pi}$.  In the sequel we use the positive square root in the expression for $p_z$, but we analyzed the other solution as well and found nothing 
different (modulo dilation and rotation).

\end{enumerate}

The first and third of these are versions of symplectic reduction.  
We thus restrict ourselves to searching for initial conditions within the embedded (non-compact) surface parametrized by $(p_x, p_y)$ with $|p_x| \leq \frac{1}{2\sqrt \pi}$.  To be clear, the constraints $x=1$ and $y=z=0$ are only imposed on the initial condition, not the entire curve.

\section{Results} \label{results}

\subsection{Main Result}

Our primary finding is the existence of a family of periodic solutions to the Kepler-Heisenberg problem, parametrized by rational numbers $j/k \in (0, 1]$, which we call the \emph{symmetry type} of the orbit.  The surprisingly beautiful trajectories are displayed in Figures \ref{k-fold} and \ref{classes}.  

\textbf{Main Numerical Discovery.} Let $j$ and $k$ be positive integers, and let $\omega=\exp (2\pi i/k)$ be the generator of the cyclic group $\mathbb{Z}_k$, acting on $\mathbb{R}^3$ by rotation about the $z$-axis by $2\pi/k$ radians.  Then we have numerically discovered a large number of solutions to the equation 
\begin{equation}  \tag{1}
q(t+T/k)=\omega^j q(t).
\end{equation}

While this result is experimental in nature, the existence of these orbits can be proved using the methods of \cite{CS}.  Indeed, the variational proof goes through essentially unchanged if one replaces the symmetry conditions (S1) and (S2) in \cite{CS} by condition (\ref{1}) above.

\subsection{Experimental approach} 

Our original goal was to find numerical approximations to the orbits whose existence was assured by Theorem {\ref{periodic}}.  
The proof of this theorem gives no clue as to the location or the initial conditions of the orbits that it asserts exist.
This motivated our choice to design a search using random initial conditions, and trying to close up the orbits which seemed quasi-periodic.  
See the first row of Figure \ref{k-fold} for examples with $k=9$ and 41.  

\begin{figure}
\centering
\begin{tabular}{cc}
\includegraphics[width=.47\textwidth]{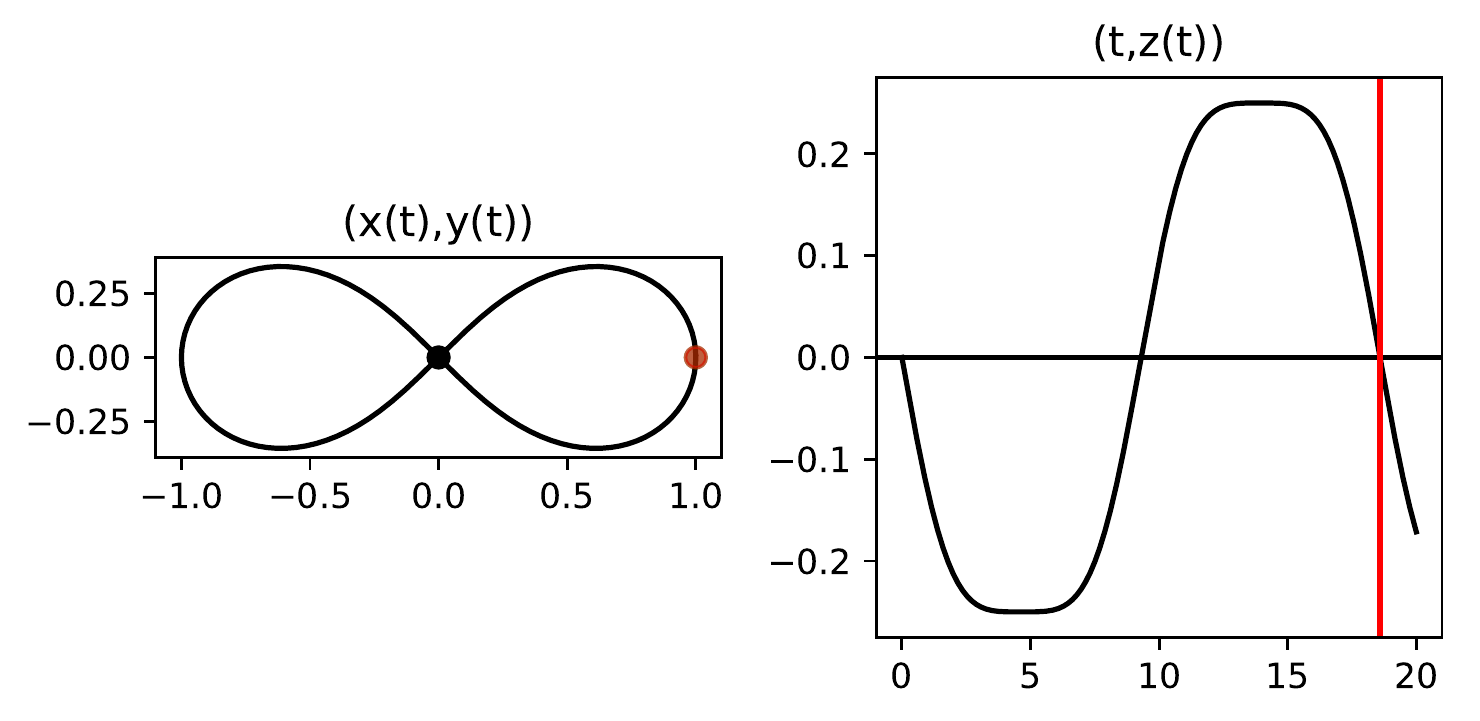} &
\includegraphics[width=.47\textwidth]{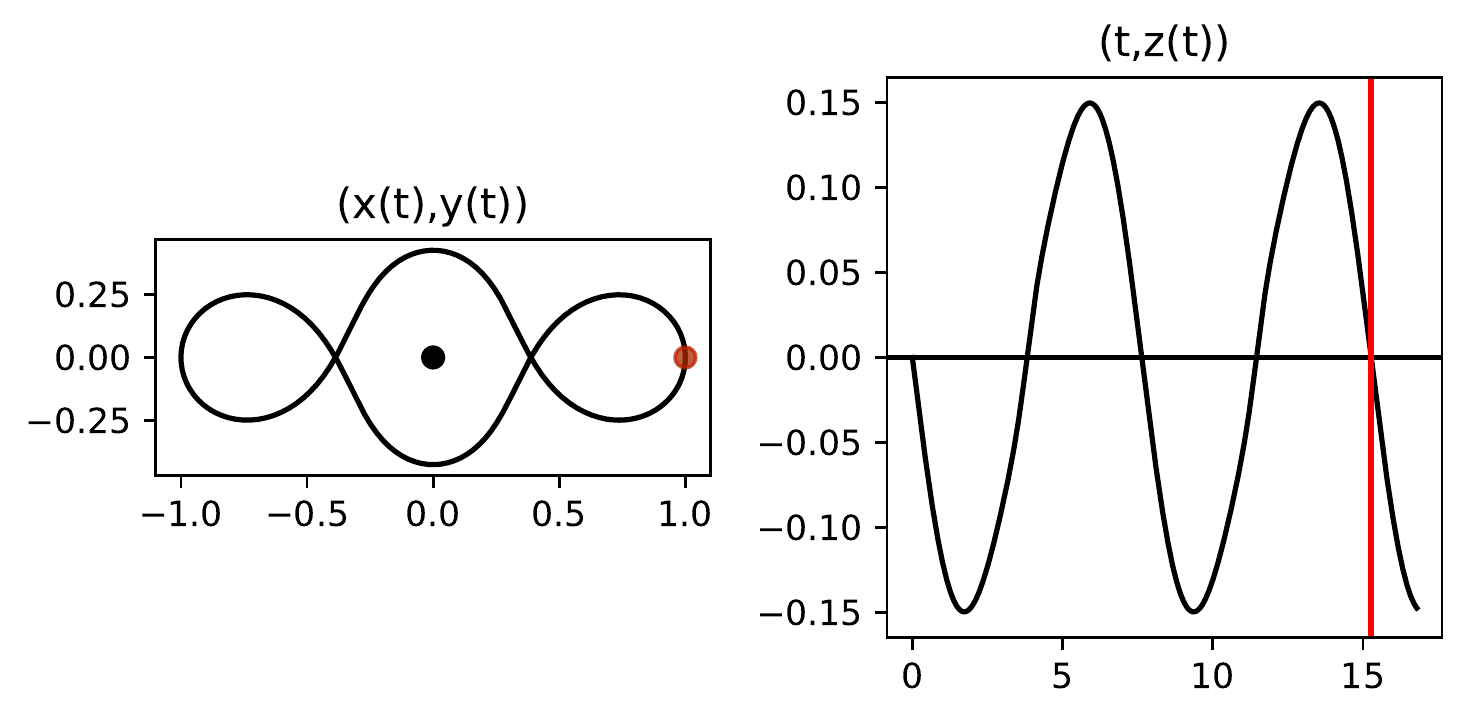}
\end{tabular}
\caption{Left: An example of a closed orbit with no rotational symmetry; the projection to the $xy$-plane exhibits 2-fold rotational symmetry, but the $z$-coordinate does not.  Right: Symmetry type $1/2$.} 
\label{1:1}
\end{figure}

\begin{figure} 
\centering
\includegraphics[width=.8\textwidth]{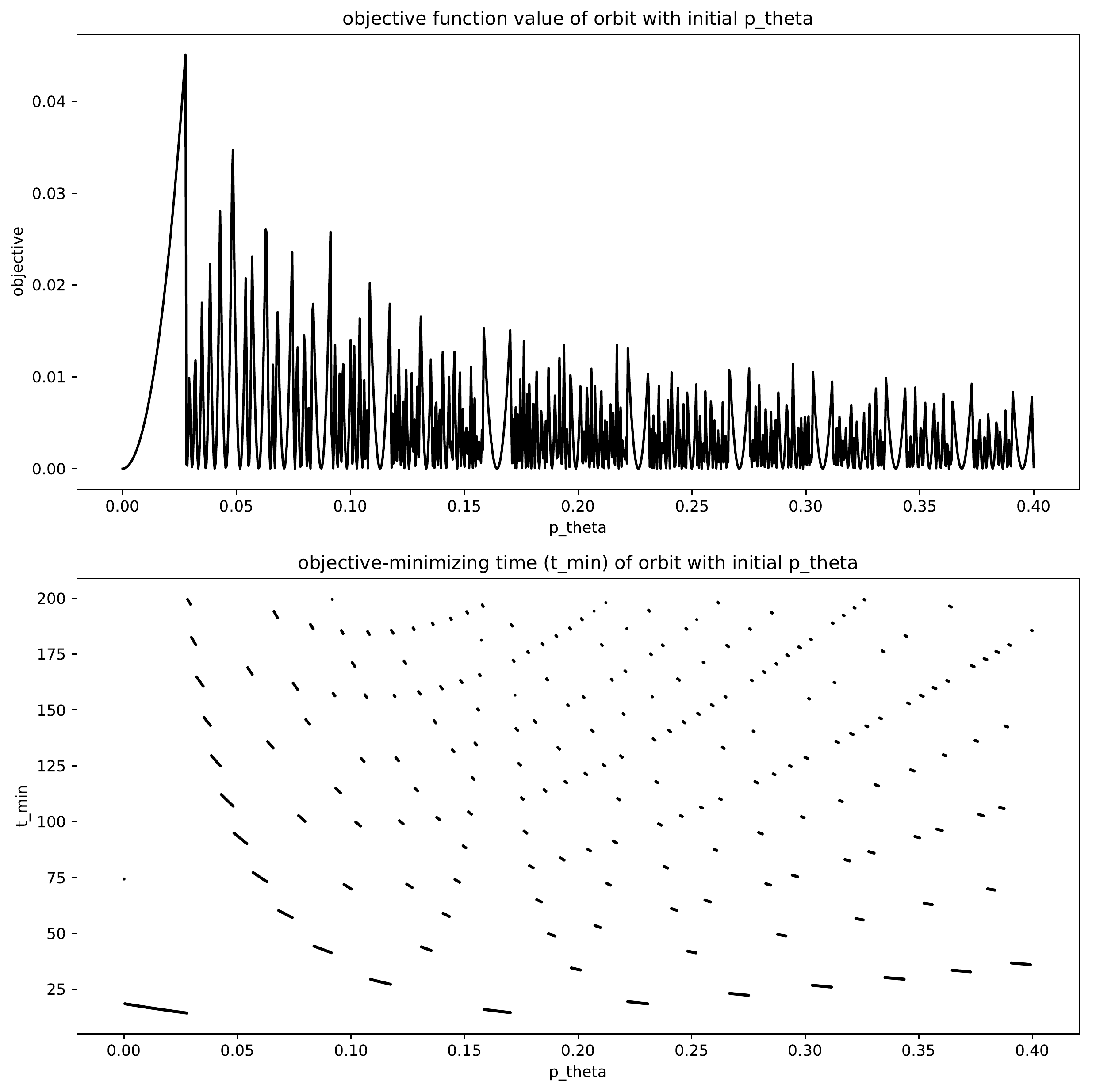}
\caption{Uniformly sampled values of $p_{\theta}\in(0,1)$ along the horizontal axes in both plots.  In the upper plot, local minima represent symmetry types of periodic orbits.  Introduced numerical thresholds obscure the fact that these are distributed densely like the rational numbers.  
In the lower plot, each discrete piece represents a symmetry type of closed orbit. Its height represents its period.  We understand how these are distributed; see Section \ref{types}.}
\label{scatter1}
\end{figure}

Our second goal was to find evidence of new closed orbits, if they exist.  
Theorem \ref{periodic} was only proved for orbits with $k$-fold symmetry for \emph{odd} $k$.  Our search yielded orbits with even degrees of symmetry as well; see the second row of Figure \ref{k-fold} for examples with $k=4$ and 8.  Moreover, we found that with a fixed $k$, there are $\varphi(k)$ many different classes of closed orbits (not just related by symmetries of the system), where $\varphi$ denotes Euler's totient function.  
Details appear below in Section \ref{types}.   We also found a closed orbit without rotational symmetry.  It appears in Figure \ref{1:1} and as an outlier in Figure \ref{scatter1}.

While we first discovered these phenomena using our \texttt{search} feature (see Section \ref{program}), which begins with randomly generated initial conditions, further investigation revealed a surprising and beautiful structure.  As described in Section \ref{reduction}, 
we reduced our search to two dimensions, which we chose to be $p_x$ and $p_y$.  Due to our choices that  $x(0)=1$ and $y(0)=z(0)=0$, we have 
$p_x(0)=J$ and $p_y(0)=p_{\theta}.$  Recall that closed orbits must satisfy $H=0$, and that this condition implies that $J$ is an integral of motion.  Therefore the two conserved quantities $J$ and $p_{\theta}$ completely determine a trajectory's motion.  

\begin{figure}[t!]
\centering
\includegraphics[width=.8\textwidth]{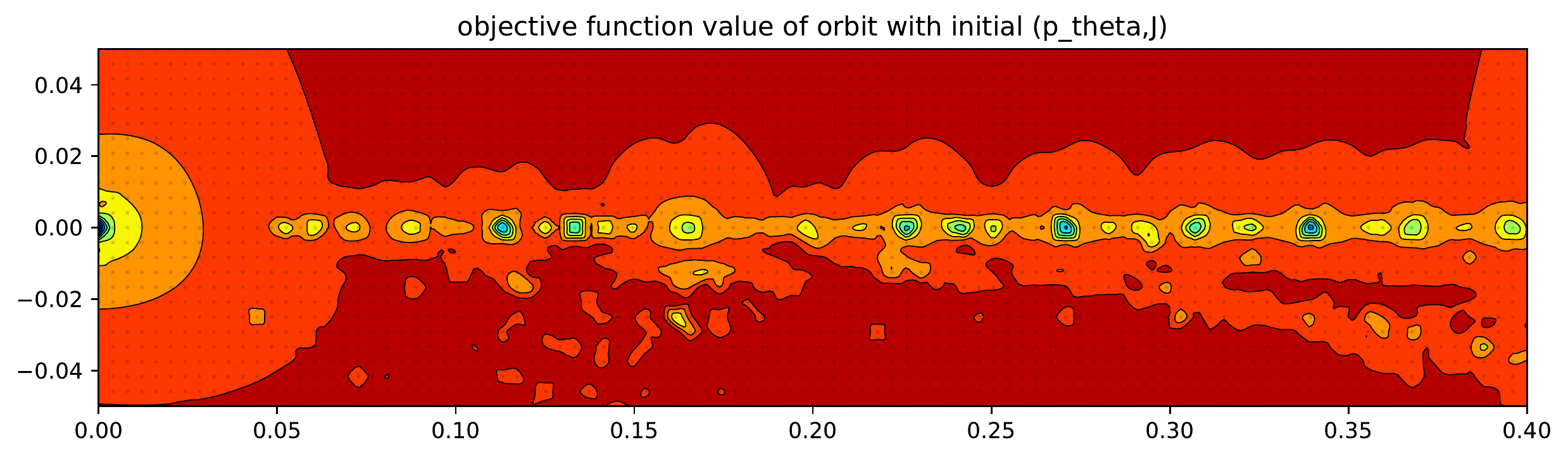}
\caption{Points in a region of the $p_{\theta}, J$-plane.  Colors represent values of the objective function, with red corresponding to large values and blue corresponding to small values.  Yellow regions lying along the line $J=0$ represent initial conditions nearby those leading to closed orbits.}
\label{scatter2}
\end{figure}

\begin{figure} 
\centering
\begin{tabular}{cc}
\includegraphics[width=.47\textwidth]{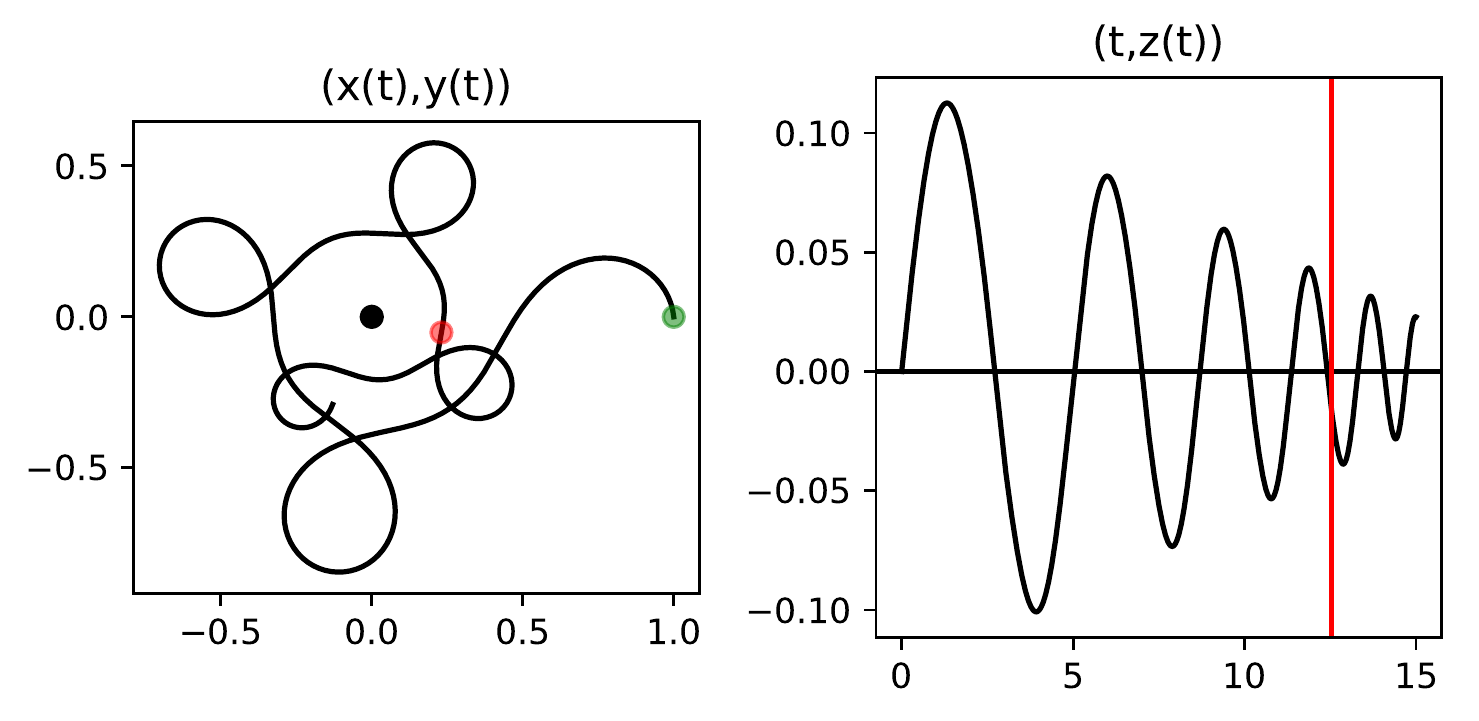} &
\includegraphics[width=.47\textwidth]{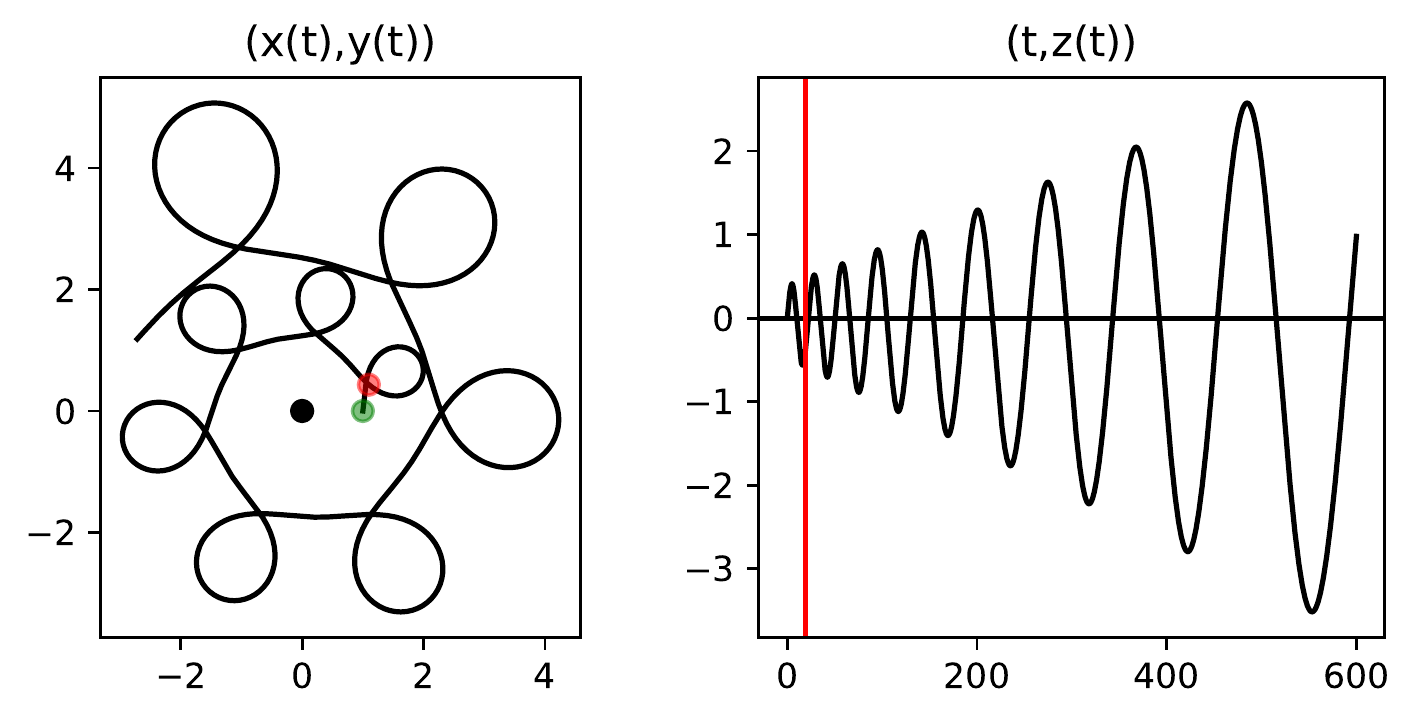}
\end{tabular}
\caption{Two examples of failed \texttt{search} procedures.  Both show orbits with non-zero dilational momentum $J$, appearing as red points in Figure \ref{scatter2}. Both orbits are self-similar in unbounded time.
Left: A quasi-periodic orbit with $J<0$ spirals in toward collision with the sun.  
Right: A potentially 8-fold periodic orbit fails to close since $J>0$, causing the orbit to dilate outwards unboundedly.  } 
\label{abortives}
\end{figure}

We uniformly sampled 2,500 points in a region of this plane and generated the scatter plot in Figure \ref{scatter2}.  The accumulation of yellow regions (signifying proximity to closed orbits) along the line $J=0$ suggested that we implement this constraint as well.
This choice can be justified as follows.  Recall that $J$ represents the \emph{dilational momentum} of an orbit, and it is conserved for closed orbits.  If $J\neq 0$, then the orbit cannot be periodic (or even bounded) -- depending on the sign of $J$, it either spirals in toward the sun or spirals out to infinity.  See Figure \ref{abortives} for examples.

Uniformly sampling 1000 points within the line $J=0$ with $p_{\theta} \in (0,1)$ generated the scatter plots in Figure \ref{scatter1}.  The lower plot shows an orbit's period as a function of $p_{\theta}$.  Each discrete piece corresponds to a particular symmetry type, described below.

\subsection{Analysis of symmetry types}
\label{types}

We intended to search for orbits with $k$-fold rotational symmetry about the $z$-axis.  After finding a large number of these, we noticed that there were multiple symmetry classes for a fixed $k$ value, which we call the \emph{order} of symmetry.  For fixed $k$, the \emph{symmetry class} is a positive integer $j<k$ which is relatively prime to $k$.  Qualitatively, if $k$ represents the number of lobes or petals in an orbit, then $j$ represents the order in which they are traced out over time.  
Also see (\ref{1}).

We call this ratio $j/k$ the \emph{symmetry type} of the orbit; it is our main invariant of interest.  We suspect that this number is a complete invariant, in the sense that this determines a closed orbit up to our rotational and dilational symmetries.  It seems nearly certain that $j/k$ is the usual rotation number from dynamical systems (see Chapter 11 of \cite{Katok} or Chapter 6 of \cite{Tabachnikov}), but we have not yet been able to formally prove this.

\begin{figure}[t!]
\centering
\begin{tabular}{cc}
\includegraphics[width=.24\textwidth]{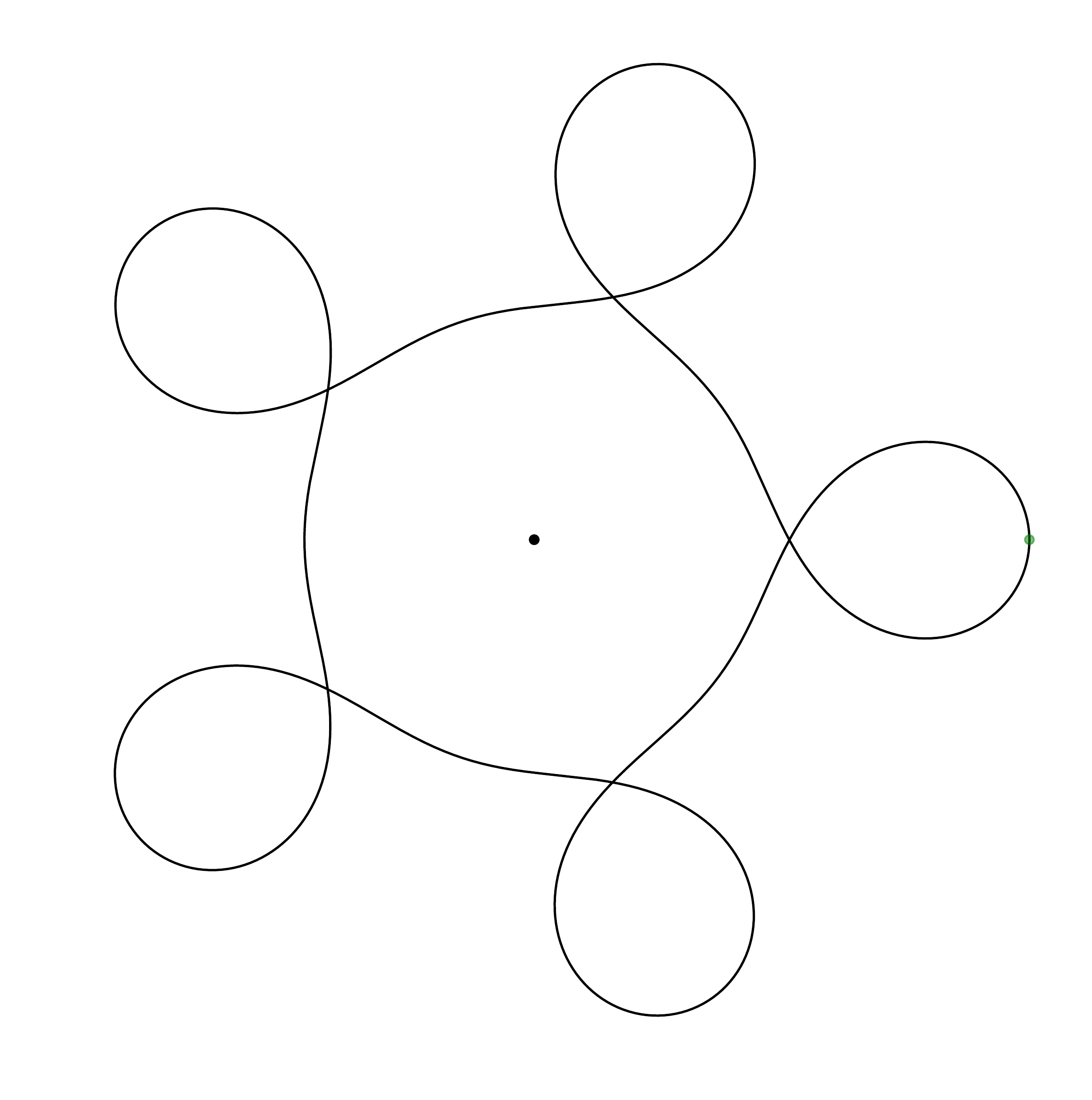}
\includegraphics[width=.24\textwidth]{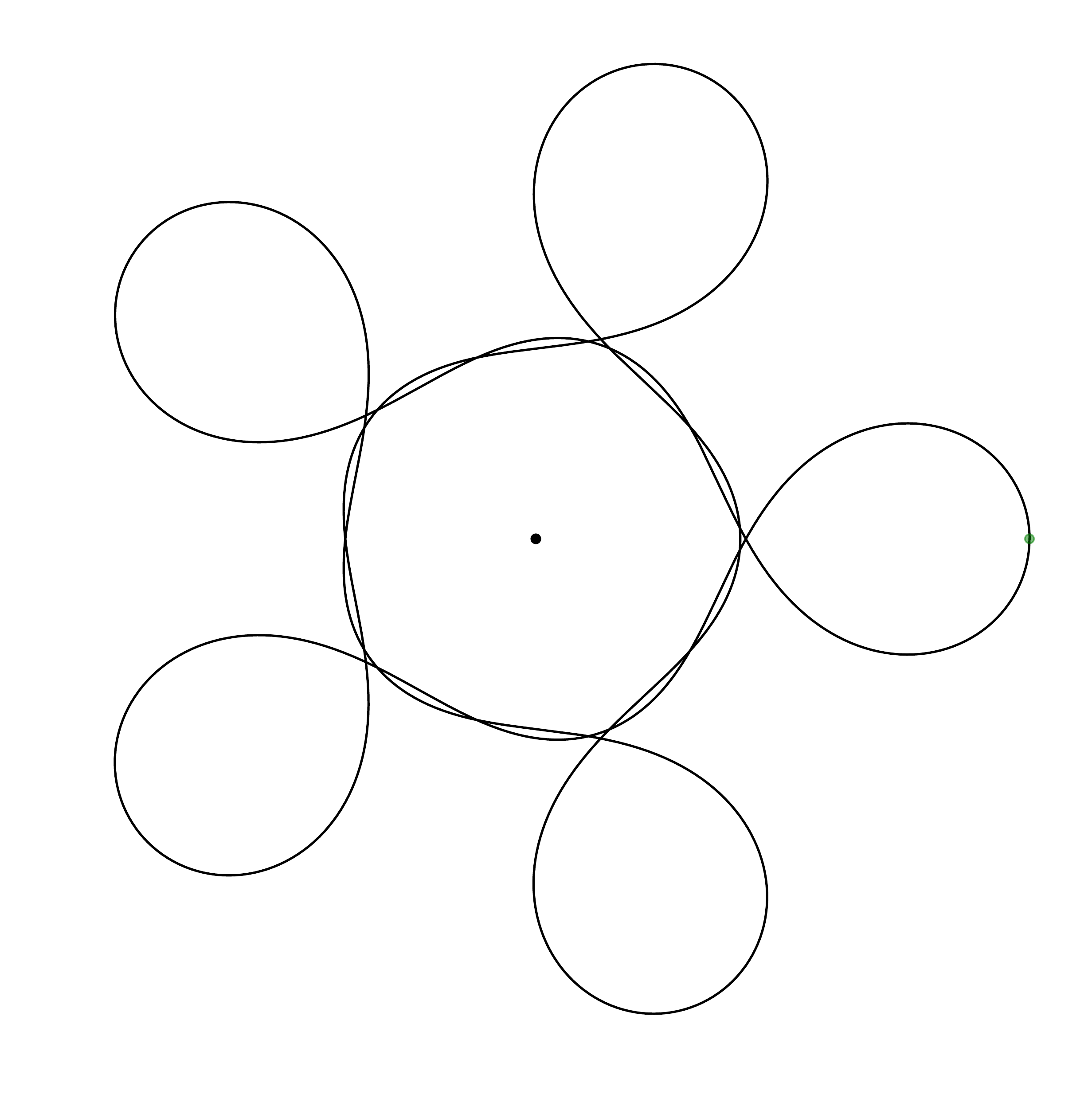}
\includegraphics[width=.24\textwidth]{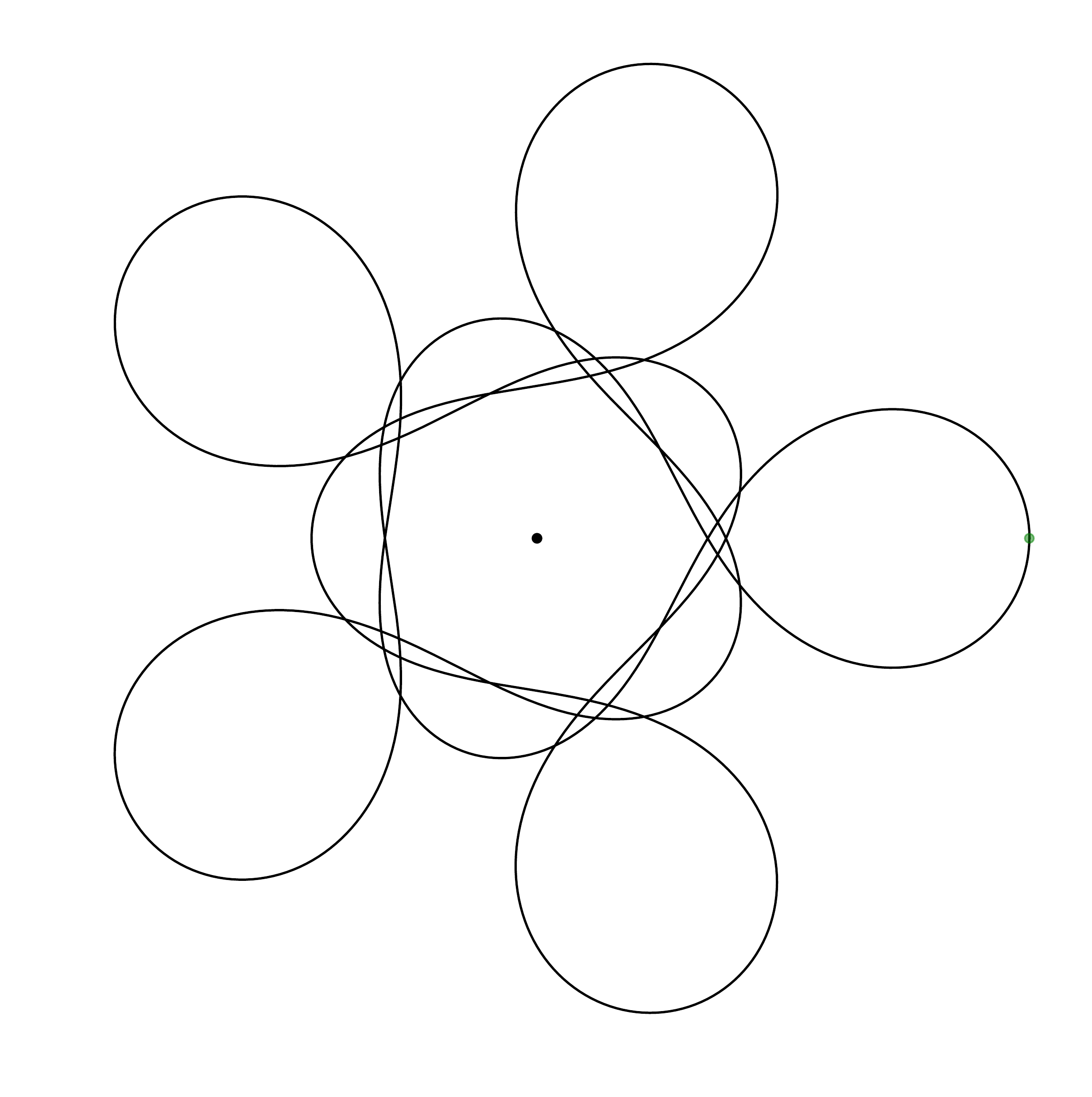}
\\
\includegraphics[width=.24\textwidth]{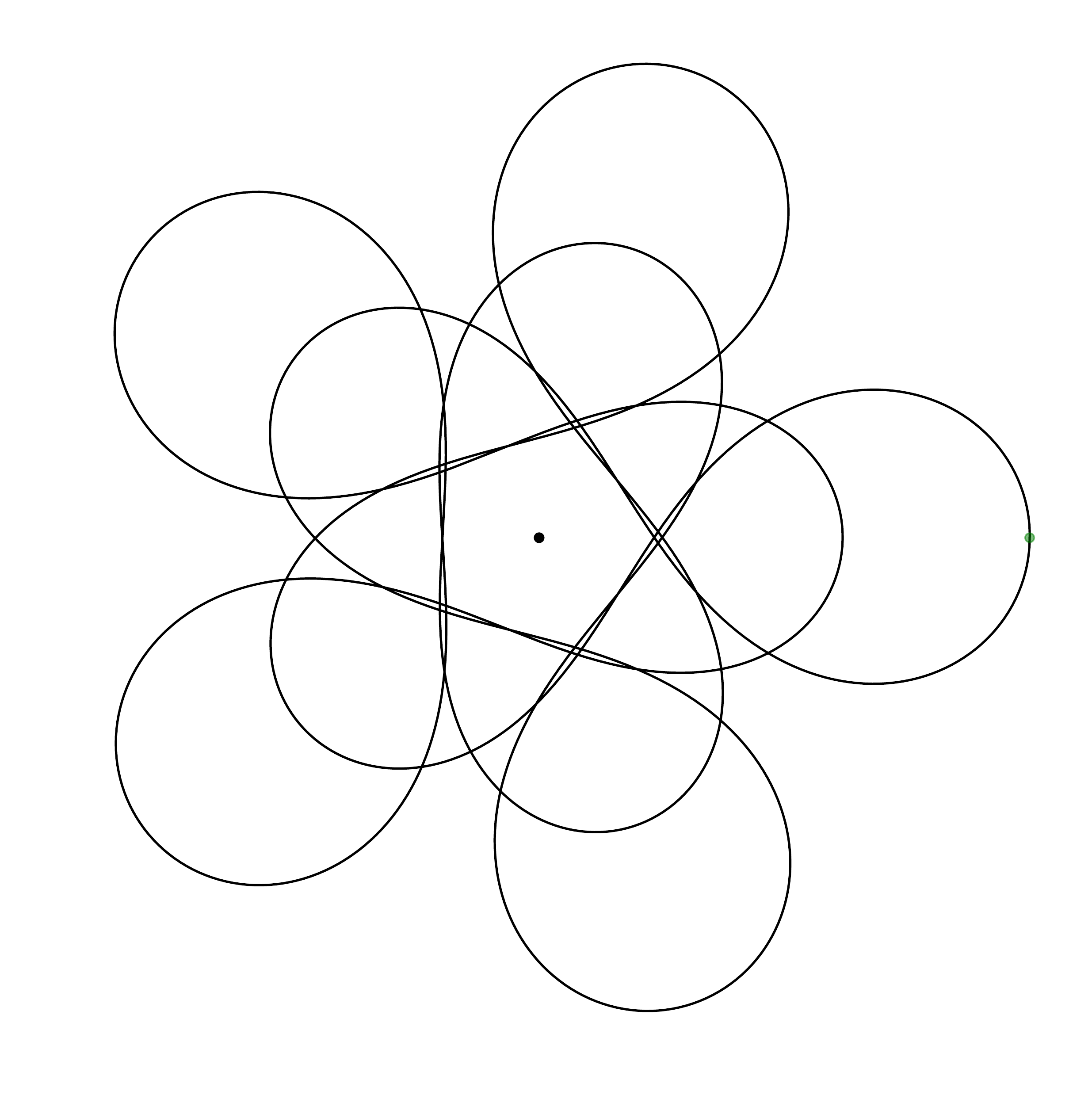}
\includegraphics[width=.24\textwidth]{order4class1.pdf}
\includegraphics[width=.24\textwidth]{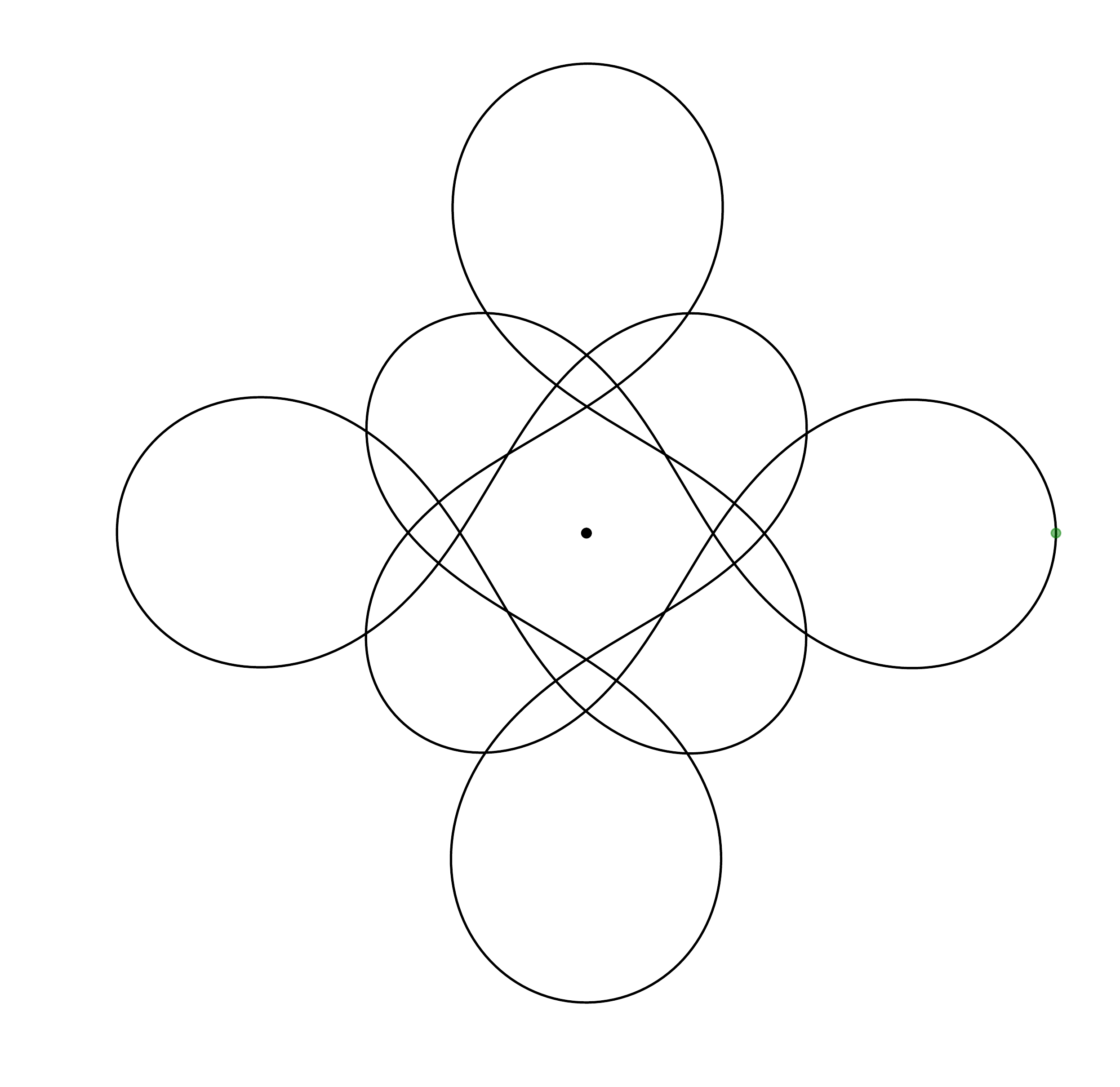}
\end{tabular}
\caption{Symmetry types $\frac{1}{5}, \frac{2}{5}, \frac{3}{5}, \frac{4}{5}, \frac{1}{4}$, and $\frac{3}{4}$.  We generated images of 175 different types.} 
\label{classes}
\end{figure}

In Figure \ref{classes}, we show all $\varphi(5)=4$ symmetry classes for $k=5$, as well as the $\varphi(4)=2$ symmetry classes for $k=4$.  Note that the orbits in Figure 
\ref{k-fold} have symmetry types $\frac{5}{9},\frac{6}{41},\frac{1}{4},\frac{7}{8}$.  
In total, we have generated images and data on over 175 different symmetry types.

There is a one-to-one correspondence between symmetry types and rational numbers in $(0,1]$.
As a function of angular momentum $p_{\theta}$, the symmetry types are arranged according to the Farey sequence; see Table \ref{table}. These symmetry types are observable in Figure \ref{scatter1}.  In the upper plot, each local minimum corresponds to a unique symmetry type.  To generate this plot we imposed a maximum period of 200, as well as other constraints, which account for the apparent finite number of minima.  In fact, the minima should densely populate the horizontal axis like $\mathbb Q \subset \mathbb R$, with the corresponding symmetry type decreasing as a function of angular momentum.

The lower plot in Figure \ref{scatter1} is most enlightening.  
Each discrete piece represents a symmetry type (individual dots are quasi-periodic orbits attracted to this symmetry type by our program, there are 1000 samples).  
The piece in the lower left corner represents symmetry type $\frac{1}{1}$; this means the orbit is periodic, but has no rotational symmetry.  See Figure \ref{1:1}.

The next lowest piece, with $p_{\theta} \approx .164$, represents the symmetry type $\frac{1}{2}$.  The ray emanating up and left from there, asymptotic to the vertical axis, consists of types $\frac{k-1}{k}$.  The ray emanating up and right from type $\frac{1}{2}$ consists of types $\frac{1}{k}$. Directly above that lies the ray of types $\frac{2}{k}$, and so on.  As all types are in reduced form, the structure in this figure is therefore somewhat number theoretic in nature.  

Our program can detect the symmetry type by using basic Fourier analysis.  The mode with highest norm of the discrete Fourier transform of the $z$-coordinate as a function of time gives the symmetry order $k$.  
We compute the discrete Fourier transform of the $xy$-curve, and the norm of the complex Fourier coefficients gives a real-valued function of the frequencies.
We periodize this function with period $k$ by taking the average over each of the congruence classes modulo $k$.  The resulting discrete function will have a maximum at $j$ if the orbit has symmetry type $\frac{j}{k}$.  Our user has the option to plot this function, the \emph{class signal}, and the discrete Fourier transforms.

\begin{table}
  \caption{
  As $p_\theta$ grows, the symmetry types are exactly the reversed Farey sequence -- in this case, the Farey sequence of order 6.
  }
\renewcommand\arraystretch{1.5}
  \begin{tabular}{ | c || c | c | c | c | c | c | c | c | c | c | c | c | }
    \hline
    $\approx p_{\theta}$ & 3.04e--6 & .060 & .071 & .087 & .113 &  .133 & .164 & .199 & .226 & .271 & .307 & .339 \\              
    \hline
    $\frac{j}{k}$ & $\frac{1}{1}$ & $\frac{5}{6}$ & $\frac{4}{5}$ & $\frac{3}{4}$ & $\frac{2}{3}$ & $\frac{3}{5}$ & $\frac{1}{2}$ & $\frac{2}{5}$ & $\frac{1}{3}$ & $\frac{1}{4}$ & $\frac{1}{5}$ & $\frac{1}{6}$ \\ \hline
  \end{tabular}
  \label{table}
\end{table}

\section{Conclusions} \label{conclusions}
\subsection{Open questions} 
Much work remains to be done.  We would like to prove that the closed orbits discovered here are indeed distributed as described in Section \ref{types}.  There may be deeper connections with rational rotation numbers and Arnold tongues.  
Also, an investigation of the knot theory of the orbits we have discovered would also be interesting, especially with respect to the symmetry types. 
And the obvious question remains: Have we found \emph{all} of the periodic orbits, or are there others outside the families described here?

For the general Kepler-Heisenberg problem, we still don't know whether the system is integrable outside of the $H=0$ subsystem.  Partial answers are given in papers like \cite{Fiorani}, which develop ``generalized" action-angle coordinates for partially integrable systems, but whether our system is completely integrable remains unknown.

\subsection{A failed attempt}

While many of our early attempts at this project ended in failure, one promising approach is worth noting despite our lack of success.
Our initial efforts were based on the methods of \cite{Nauenberg1, Nauenberg2}, which use techniques from Fourier analysis.  We attempted to take advantage of the variational formulation and work within the Lagrangian formalism.  
The sub-Riemannian structure on the Heisenberg group allows one to parameterize the optimization problem in terms of a single complex-valued curve, the Fourier decomposition of which lends itself to a particularly simple expression of the symmetry conditions.  
Our approach was roughly as follows.

\begin{enumerate}
\item Input Fourier coefficients with symmetry for plane curve $(x,y)$.
\item Construct $z$ = periodic + linear. Let $Q$ be proportional to the linear part of $z$. Want $Q=0$ so $z$ periodic. 
\item Reconstruct parametrized space curve. 
\item Compute action of path. 
\item Let $A$ = action + $\lambda Q$ be the constrained action. Want critical point of $A$. 
\item Minimize norm of $\nabla A$ as function of Fourier coefficients and $\lambda$. 
\item Output solution to equations of motion with prescribed $(x,y)$ symmetry and periodic $z$. 
\end{enumerate}
While this approach was conceptually attractive, it was surprisingly difficult to implement.  In \cite{Nauenberg1, Nauenberg2}, configuration space is 2-dimensional.  The complications added by our third dimension -- large numbers of composed functions and constrained optimization in relatively high dimension -- led us to abandon this idea for the simpler approach described above.  The interested reader can find the discarded code in \cite{Github}.

\subsection{Codebase}


All code used for the experiments and plots in this work is hosted publicly as free, open-source software at \url{https://github.com/vdods/heisenberg}. Explicit instructions to reproduce each plot in this paper are provided there, as are all relevant documentation and licensing information.  This code-publishing is an  explicit effort to conduct experimental mathematics so that all results are 100\% available and exactly reproducible.

\subsection*{Acknowledgements}

The authors are very grateful to the referee for suggestions which significantly improved the exposition and focus of this paper.

\newpage

\end{document}